\documentclass[a4paper,12pt,reqno]{amsart}
\usepackage[a4paper,margin=1.9cm,top=2.5cm,bottom=2.5cm,centering,vcentering]{geometry}
\usepackage{hyperref,graphicx,tikz,amsmath,amssymb,amsfonts,amsthm}
\usepackage{makecell}
\setcellgapes{5pt}

\definecolor{webgreen}{rgb}{0,.5,0}
\definecolor{webbrown}{rgb}{.6,0,0}
\definecolor{red}{rgb}{1,0,0}

\usepackage[normalem]{ulem}

\usepackage{color}
\usepackage{url}

\usetikzlibrary{decorations.pathreplacing}

\numberwithin{equation}{section}
\theoremstyle{plain}
\newtheorem{theorem}{Theorem}[section]

\newtheorem{definition}[theorem]{Definition}

\newtheorem{corollary}[theorem]{Corollary}
\newtheorem{proposition}[theorem]{Proposition}

\begin{document}
	
\title{Partial sums of the Gibonacci sequence}
   
\author[P. J. Mahanta]{Pankaj Jyoti Mahanta}
\address{Gonit Sora, Dhalpur, Assam 784165, India}
\email{pankaj@gonitsora.com}


\begin{abstract}
	Recently, Chu studied some properties of the partial sums of the sequence $P^k(F_n)$, where $P(F_n)=\big(\sum_{i=1}^nF_i\big)_{n\geq1}$ and $(F_n)_{n\geq1}$ is the Fibonacci sequence, and gave its combinatorial interpretation. We generalize those results, introduce colored Schreier sets, and give another equivalent combinatorial interpretation by means of lattice path.
\end{abstract}

\maketitle

\section{Introduction}
The Fibonacci sequence is defined by $F_n=F_{n-1}+F_{n-2}$, with initial terms $F_1=1$ and $F_2=1$. One of its generalizations is the Gibonacci sequence, which is defined by $G_n=G_{n-1}+G_{n-2}$, with initial terms $G_1$ and $G_2$, where the initial two terms can be any positive integer. Both the sequences are related by the following identity
\[G_n=G_1F_{n-2}+G_2F_{n-1} \ \text{for all} \ n>2.\]
We refer the reader to Benjamin and Quinn’s book \cite{benjamin2003proofs}, and the author’s joint work with Saikia \cite{msGibo2021}.

In a recent paper \cite{chu2021partial}, Chu defined a function $P$ such that
\[P(F_n):=\bigg(\sum_{i=1}^nF_i\bigg)_{n\geq1}.\]
Chu also defined that
\[P^k(F_n):=P(P^{k-1}(F_n)) \ \text{for all} \ k\geq2,\]
and denoted the $n$th term of the sequence $P^k(F_n)$ by $a_k(n)$. We generalize these to Gibonacci sequence, that is
\[P(G_n):=\bigg(\sum_{i=1}^nG_i\bigg)_{n\geq1} \ \text{and} \ P^k(G_n):=P(P^{k-1}(G_n)) \ \text{for all} \ k\geq2.\]
We denote the $n$th term of the sequence $P^k(G_n)$ by $a_k^{\{G_1,G_2\}}(n)$. So, $a_k^{\{F_1,F_2\}}(n)=a_k(n)$. For simplicity, sometimes we write $a_k^\prime(n)$ instead of $a_k^{\{G_1,G_2\}}(n)$. So, $a_k^\prime(n)=a_k^\prime(n-1)+a_{k-1}^\prime(n)$.

\section{Generalization of $a_k(n)$}
We start with the following interesting binomial coefficient identity.
\begin{proposition}\label{pro1}
For all non-negative integers $k$ and $n$,
\[\sum_{i=0}^{n}\binom{k+i}{k}=\binom{k+n+1}{k+1}.\]
\end{proposition}
It can be easily proved by mathematical induction on $n$ using the Pascal's identity $\dbinom{n}{k}+\dbinom{n}{k+1}=\dbinom{n+1}{k+1}$. Note that $\dbinom{0}{0}=1$, and if $n\geq0$ and $\ell<0$ then $\dbinom{n}{\ell}=0$.

\begin{theorem}
For all integers $n,k\geq1$, we have
\[a_k^{\{G_1,G_2\}}(n)= \sum_{i=0}^{n-1}\binom{k-1+i}{k-1}G_{n-i}.\]
\end{theorem}

The  theorem can be easily proved by mathematical induction using Proposition \ref{pro1}.

\begin{corollary}[Generalization of Lemma 2.1, \cite{chu2021partial}]\label{CorT1}
For $k\geq0$, we have
\[a_k^\prime(3)=(k+1)G_2+\bigg(\binom{k+1}{k-1}+1\bigg)G_1.\]
\end{corollary}

\begin{corollary}[Generalization of Theorem 1.1, \cite{chu2021partial}]\label{CorT2}
For all $n,k\geq1$, we have
\[a_k^\prime(n)=a_{k-1}^\prime(n+2)-\binom{n+k-1}{k-1}G_2-\binom{n+k-1}{k-2}G_1.\]
\end{corollary}

When $k=1$, it gives us the well-known identity $\displaystyle\sum_{i=1}^nG_i=G_{n+2}-G_2$, and then $\displaystyle\sum_{i=1}^nF_i=F_{n+2}-1$.

\begin{proof}[Proof of Corollary \ref{CorT2}]
We have,
\begin{align*}
a_{k-1}^\prime(n+2) & =\sum_{i=0}^{n+1}\binom{k-2+i}{k-2}G_{n+2-i}\\
& =G_{n+2} +(k-1)G_{n+1}+ \sum_{i=0}^{n-1}\binom{k+i}{k-2}G_{n-i}\\
& = kG_2 + \bigg(1+\binom{k}{k-2}\bigg)G_n+ \sum_{i=1}^{n-1}\bigg(k+\binom{k+i}{k-2}\bigg)G_{n-i}.\\
&\hspace{20pt} \big(\text{Since we get} \ G_{n+2}=\sum_{i=0}^{n-1}G_{n-i}+G_2 \ \text{and} \ G_{n+1}=\sum_{i=1}^{n-1}G_{n-i}+G_2.\big)
\end{align*}

Therefore, $a_{k-1}^\prime(n+2)-a_k^\prime(n)$ is equal to

\begin{align*}
& kG_2 + \binom{k}{k-2}G_n + \sum_{i=1}^{n-1}\bigg(k+\binom{k+i}{k-2}-\binom{k-1+i}{k-1}\bigg)G_{n-i}\\
=& \bigg(k+\binom{k}{k-2}\bigg)G_2 + \bigg(k+\binom{k+1}{k-2}-\binom{k}{k-1}\bigg)G_{n-1}\\ &\hspace{140pt}+\sum_{i=2}^{n-1}\bigg(k+\binom{k}{k-2}+\binom{k+i}{k-2}-\binom{k-1+i}{k-1}\bigg)G_{n-i}\\
&\hspace{200pt} \big(\text{Since} \ G_n=\sum_{i=2}^{n-1}G_{n-i}+G_2.\big)\\
=& \binom{k+1}{k-1}G_2 + \binom{k+1}{k-2}G_{n-1} + \sum_{i=2}^{n-1}\bigg(\binom{k+1}{k-1}+\binom{k+i}{k-2}-\binom{k-1+i}{k-1}\bigg)G_{n-i}\\
=& \binom{k+2}{k-1}G_2 + \binom{k+2}{k-2}G_{n-2} + \sum_{i=3}^{n-1}\bigg(\binom{k+2}{k-1}+\binom{k+i}{k-2}-\binom{k-1+i}{k-1}\bigg)G_{n-i}\\
&\hspace{200pt} \big(\text{Since} \ G_{n-1}=\sum_{i=3}^{n-1}G_{n-i}+G_2.\big)\\
=& \binom{k+3}{k-1}G_2 + \binom{k+3}{k-2}G_{n-3} + \sum_{i=4}^{n-1}\bigg(\binom{k+3}{k-1}+\binom{k+i}{k-2}-\binom{k-1+i}{k-1}\bigg)G_{n-i}.\\
&\hspace{200pt} \big(\text{Since} \ G_{n-2}=\sum_{i=4}^{n-1}G_{n-i}+G_2.\big)
\end{align*}
Proceeding in this way upto $(n-2)$ steps we get that the difference is equal to

\[\binom{k+n-2}{k-1}G_2 + \binom{k+n-2}{k-2}G_{2} + \bigg(\binom{k+n-2}{k-1}+\binom{k+n-1}{k-2}-\binom{k-1+n-1}{k-1}\bigg)G_{1},\]
which is equal to
\[\binom{n+k-1}{k-1}G_2+\binom{n+k-1}{k-2}G_1.\]
\end{proof}

\section{A combinatorial interpretation of $a_k^{\{G_1,G_2\}}(n)$}\label{Schreier}
A finite subset $S$ of natural numbers is called a Schreier set if $\min S\geq |S|$. By counting some Schreier sets Chu gave a combinatorial interpretation of $a_k(n)$. We generalize it for $a_k^{\{G_1,G_2\}}(n)$.

\begin{definition}[Chu]
For any integers $n\geq 1$, and $k\geq 0$,
\[s_k(n):=\#\{S\subset \{1,2,3,\dots,n\}: |S|\geq k, \ \text{and} \ \min S\geq |S|\}.\]
\end{definition}
It is easy to prove the following proposition.
\begin{proposition} For any natural number $n$,
\[\#\{S\subset \{1,2,3,\dots,n\}: |S|=\ell, \ \text{and} \ \min S>\ell\}=\binom{n-\ell}{\ell},\]
and
\[\#\{S\subset \{1,2,3,\dots,n\}: |S|=\ell, \ \text{and} \ \min S=\ell\}=\binom{n-\ell}{\ell-1}.\]
\end{proposition}

This implies that
\[s_k(n)=\sum_{\ell\geq k}\bigg(\binom{n-\ell}{\ell}+\binom{n-\ell}{\ell-1}\bigg).\]

Now we count some Schreier sets where in each set, one specific element is of different type. We call these sets colored Schreier sets. Here, for the sets in \[\{S\subset \{1,2,3,\dots,n\}: |S|=\ell, \ \text{and} \ \min S=\ell\},\] the element $\ell$ occurs in $G_1$ different colors. And for the sets in \[\{S\subset \{1,2,3,\dots,n\}: |S|=\ell, \ \text{and} \ \min S>\ell\},\] the element which is equal to $|S|$ occurs in $G_2$ different colors. For example, for $G_1=3$ and $G_2=2$, let $^R$, $^B,$ and $^G$ be three different colors, and then all the colored Schreier sets corresponding to $\{S\subset \{1,2,3,4,5,6\}: |S|\geq 2, \ \text{and} \ \min S\geq |S|\}$ are

\begin{flushleft}
$\{2^R,3\},\{2^R,4\},\{2^R,5\},\{2^R,6\},\hspace{45pt} \{3^R,4\},\{3^R,5\},\{3^R,6\},\{4^R,5\},\{4^R,6\},\{5^R,6\},$

$\{2^B,3\},\{2^B,4\},\{2^B,5\},\{2^B,6\},\hspace{44pt} \{3^B,4\},\{3^B,5\},\{3^B,6\},\{4^B,5\},\{4^B,6\},\{5^B,6\},$

$\{2^G,3\},\{2^G,4\},\{2^G,5\},\{2^G,6\},$

$\{3^R,4,5\},\{3^R,4,6\},\{3^R,5,6\},\hspace{53pt} \{4^R,5,6\},$

$\{3^B,4,5\},\{3^B,4,6\},\{3^B,5,6\},\hspace{52pt} \{4^B,5,6\},$

$\{3^G,4,5\},\{3^G,4,6\},\{3^G,5,6\}.$
\end{flushleft}
The total number of colored Schreier sets for particular $n$ and $k$ is given by
\[s_k^{\{G_1,G_2\}}(n):=\sum_{\ell\geq k}\bigg(\binom{n-\ell}{\ell}G_2+\binom{n-\ell}{\ell-1}G_1\bigg).\]
For simplicity, sometimes we denote it by $s_k^\prime(n)$.

\begin{proposition} [Generalization of Corollary 2.3, \cite{chu2021partial}]
For $k\geq 0$ and $n\geq 1$, we get
\[s_{k+1}^{\{G_1,G_2\}}(n)=s_k^{\{G_1,G_2\}}(n)-\binom{n-k}{k}G_2-\binom{n-k}{k-1}G_1.\]
\end{proposition}

\begin{theorem} [Generalization of Theorem 1.3, \cite{chu2021partial}]
For $k\geq0$ and $n\geq 1$, we get
\[s_k^{\{G_1,G_2\}}(n)=a_k^{\{G_1,G_2\}}(n-2(k-1)).\]
\end{theorem}

\begin{proof} First we show that $s_0^\prime(n)=a_0^\prime(n+2)=G_{n+2}$. We have,
\begin{align*}
s_0^\prime(n) &=\sum_{\ell\geq 0}\bigg(\binom{n-\ell}{\ell}G_2+\binom{n-\ell}{\ell-1}G_1\bigg)\\
&=F_{n+1}G_2+F_nG_1\\
&\hspace{20pt} (\text{Since for} \ n\geq1, F_n=\sum_{i=0}^{\lfloor\frac{n-1}{2}\rfloor}\binom{n-i-1}{i}, \ \text{which is a well-known identity}.)\\
&=G_{n+2}.
\end{align*}

The remaining part of the proof is similar to that of the Theorem 1.3 of \cite{chu2021partial}.
\end{proof}

\section{Combinatorial interpretation of $a_k^{\{G_1,G_2\}}(n)$ by means of lattice paths}\label{latice}
\begin{figure}[htb!]
	\includegraphics[width=0.99\textwidth]{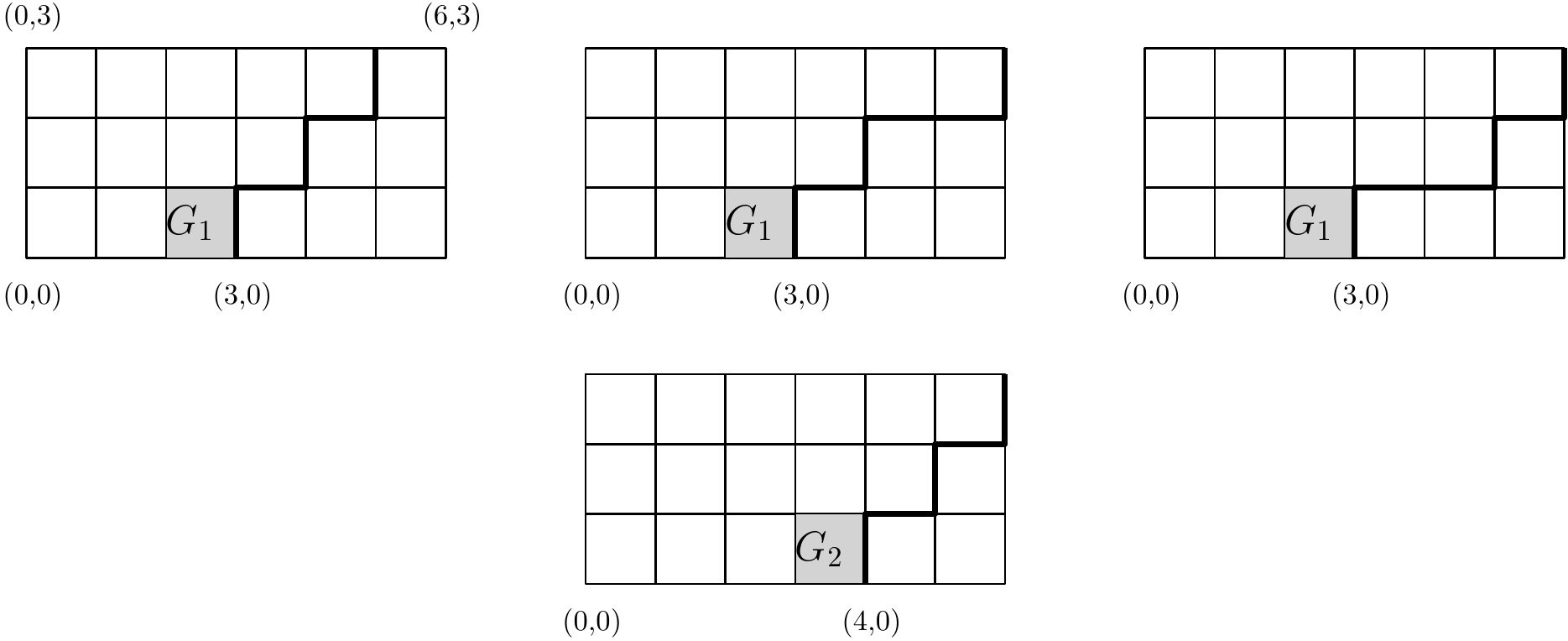}
	\caption{Lattice paths for $s_3^\prime(6)=a_3^\prime(2)$.} \label{fig1}
\end{figure}

Let us define a set of lattice paths in a $k\times n$ grid, which start from $(\ell,0)$, where $\ell\geq k$, and end on the line joining $(0,k)$ and $(n,k)$, and which consist only of steps in the upward or rightward directions, such that
\begin{itemize}
	\item the first step is always in upward direction,
\item only one step can be taken at once in the upward direction,
\item one or more steps can be taken at once in the rightward direction.
\end{itemize}
The square located in the $i$th column and the $j$th row from the lower left corner of a grid is called the $(i,j)$-cell of the grid. If $\ell=k$ then we give $G_1$ colors to the $(k,1)$-cell, if $\ell>k$ then we give $G_2$ colors to the $(\ell,1)$-cell, and when $k=0$ then we consider that only one lattice path is there along with $G_2$ colors. We observe that, the total lattice paths for $0\leq k\leq \bigg\lfloor\dfrac{n+1}{2}\bigg\rfloor$ is equal to $s_k^{\{G_1,G_2\}}(n)=a_k^{\{G_1,G_2\}}(n-2(k-1))$. For example, Figure \ref{fig1} shows all lattice paths for $s_3^\prime(6)$.

If a lattice path starts from $(0,0)$ and ends at $(n,k)$ in a $k\times n$ grid, and it consists only steps in the upward or rightward directions, then the total number of such lattice paths is equal to $\dbinom{n+k}{k}$. We can construct all such lattice paths for $0\leq k\leq \bigg\lfloor\dfrac{n+1}{2}\bigg\rfloor$. By Lemma 2.2 of \cite{chu2021partial} we get, if $G_1=F_1$ and $G_2=F_2$, then the set of the above lattice paths is in one-to-one correspondence with the set of these lattice paths in a $k\times (n-2k+1)$ grid for $0\leq k\leq \bigg\lfloor\dfrac{n+1}{2}\bigg\rfloor$. For example, the total lattice paths of this type in a $3\times1$ grid is 4, and from Figure \ref{fig1} we get $s_3^{\{F_1,F_2\}}(6)=4$.

\section*{Acknowledgements}
The author would like to thank Manjil P. Saikia for his helpful comments.

\bibliographystyle{alpha}
\bibliography{ref}

\end{document}